\documentstyle{amsart}
\newcommand{\bext}{\mbox{$Bext^2(G,T)$}}
\newtheorem{theorem}{Theorem}
\newtheorem{lemma}[theorem]{Lemma}
\newtheorem{claim}[theorem]{Claim}
\newtheorem{cor}[theorem]{Corollary}
\newtheorem{definition}{Definition}
\newcommand{\boldt}{\mbox{\boldmath $t$}}
\newcommand{\pnmu}{\mbox{$p^n_\mu$}}
\newcommand{\boldQ}{\mbox{\boldmath $Q$}}

\begin{document}
\title[\bext\ UNDER GCH]{\bext\ can be Nontrivial, even assuming GCH}
\author[M. Magidor]{Menachem Magidor}
\address{Institute of Mathematics\\Hebrew University, Jerusalem Israel 91904}
\email{menachem@@math.huji.ac.il}
\author[S. Shelah]{Saharon Shelah}
\address{Institue of Mathematics\\Hebrew University, Jerusalem Israel 91904}
\email{shelah@@math.huji.ac.il}
\thanks{The research of the second named author was supported by the Basic
Research Foundation of the Israeli Academy of Science (Publ. No. 514).}
\thanks{This paper is in final form and no version of it will be submitted
for publication elsewhere.}
\subjclass{20K20, 20K40;Secondary 03E55, 03E35}
\keywords{Butler Groups, Balanced Extensions, Consistency Proofs}
\bibliographystyle{amsplain}
\maketitle
\begin{abstract} Using the consistency of some large cardinals we produce
a model of Set Theory in which the generalized continuum hypothesis
holds and for some torsion-free abelian group $G$ of cardinality
$\aleph_{\omega+1}$ and
for some torsion group $T$
$$\bext\not=0.$$
Hence G.C.H. is not sufficient for getting the results of
\cite{Fu-Ma}.
\end{abstract}
\section{Introduction}
All groups in this paper are abelian groups. For basic terminology about
abelian groups in general we refer the reader to \cite{Fu.book}. For
terminology concerning Butler groups see \cite{Bi-Sa,Al-Hi,DHR,Fu-Ma,Fu-new}.
It is commonly agreed that the three major questions concerning the infinite
rank Butler groups are:
\begin{enumerate}
\item Are $B_1$-groups necessarily $B_2$-groups?
\item Does $\bext=0$ hold for all torsion-free groups $G$ and torsion
groups $T$?
\item Which pure subgroups of $B_2$-groups are again $B_2$-groups? In
particular: is a balanced subgroup of a $B_2$-group a $B_2$-group?
\end{enumerate}

In \cite{Bi-Sa} it is shown that the answer to all these questions is
``Yes" for countable groups $G$. In the series of papers \cite{Al-Hi,DR,DHR}
it was shown that under the continuum hypothesis the answer is ``Yes" to
all three questions for groups $G$ of cardinality $\le\aleph_\omega$. In
\cite{Du-Th} it is shown that the answer to question 2 is ``No" if the
continuum hypothesis fails. In a more recent paper \cite{Fu-Ma} it is
shown that in the constructible universe, $L$ the answer is ``Yes" to all
three questions for {\em arbitrary} groups $G$. Actually \cite{Fu-Ma}
used only the generalized continuum hypothesis and that the
combinatorial principle $\Box_\kappa$ holds for every singular cardinal
$\kappa$ whose cofinality is $\aleph_0$. Is the use made in \cite{Fu-Ma}
of the additional combinatorial principle really needed or does the
affirmitive answer to our three questions follow simply from G.C.H.?
Let us mention that a key tool used in \cite{DHR,Fu-Ma} was the
representation of an arbitrary torsion-free group as the union of a chain of
subgroups which are countable unions of balanced subgroups. In \cite{For-Ma}
it is shown that such a representation is equivalent to a weak version of
$\Box_\kappa$.

In this paper we show that at least for getting an affirmitive answer to
questions 2 and 3, one needs some extra set theoretic assumptions in addition
to
G.C.H.\@ We do it by producing a model of Set Theory, satisfying G.C.H., in
which for some torsion-free $G$ of cardinality $\aleph_{\omega+1}$ and some
torsion $T$, $\bext\not=0$. Also in the same model there will be a balanced
subgroup of a completely decomposable group which is not a $B_2$-group. Hence
the answer to question 3 in this model is ``No".
The construction of the model requires the
consistency of some large cardinals, which can not be avoided since getting a
model in which $\Box_\kappa$ fails for some singular $\kappa$ requires
assumptions stronger than the consistency of Set Theory. Let us stress that
the status of question 1 is not known and it is possible (though unlikely)
that the implication ``every $B_1$-group is a $B_2$-group" is a theorem of
Set Theory.

Since this paper is aimed at a mixed audience of set theorists and abelian
group theorists it is divided into two sections with very different
prerequisites. In the next section we
describe the construction of the model of Set Theory with certain properties
to be listed below. In the following section we shall describe how to use
the listed properties to get a group $G$ which will be the counterexample
to $\bext=0$. A reader who is not familiar with standard set theoretical
techniques, like forcing, can skip the set theoretic section and simply
assume the properties of the model listed below. We do assume some basic
Set Theory at the level introduced by \cite{Ek-Me}.

We now describe the properties of the model which will be used in the
construction of the counterexample to questions 2 and 3. The model will
naturally satisfy G.C.H.\@ Hence by standard cardinal arithmetic
$\aleph_{\omega}^{\aleph_0}=\aleph_{\omega+1}$. Therefore we can enumerate
all the $\omega$-sequences from $\aleph_\omega$ in a sequence of order type
$\aleph_{\omega+1}$. Let $\langle f_\alpha|\alpha<\aleph_{\omega+1}\rangle$
be this enumeration. Let $F_\alpha$ be the range of $f_\alpha$.
The important property of the model is the following:

{\em For some stationary subset $S$ of $\aleph_{\omega+1}$ such that every
point of $S$ has cofinality $\aleph_1$,and for some choice of a cofinal
set $C_\beta$ in $\beta$ of order type $\omega_1$, for every $\beta\in S$ and
for some fixed countable
ordinal $\delta$ we have:
\begin{enumerate}
\item $$\bigcup_{\alpha\in D}F_\alpha$$ has order type $\delta$ for every
$D\subseteq C_\beta$ which is
cofinal subset of $C_\beta$ and for every $\beta\in S$. In particular for
$D=C_\beta$ $$E_\beta=\bigcup_{\alpha\in C_\beta}F_\alpha$$ has order type
$\delta$.
\item If $\beta\not=\gamma$ both in $S$, then $E_\beta\cap E_\gamma$ has
order type less than $\delta$.
\item $\delta$ is an indecomposable ordinal, namely $\delta$ can not be
represented
as a finite sum of smaller ordinals. Or equivalently, $\delta$ is not the
finite union of sets of ordinals of order type less than $\delta$.
\end{enumerate}
}

Denote the conjunction of all the properties above by (*).
The main theorem of Section 1 is
\begin{theorem}\label{main1} Assume the consistency of a supercompact
cardinal. Then
there is a model of Set Theory in which (*) holds. The model also satisfies
the Generalized Continuum Hypothesis.
\end{theorem}
The construction of the model is very close to the construction in
\cite{HJSh249}.
The main tool that will be used to get in Section 3 an example of a group
G
satisfying $ \bext\not=0$ is the notion of $\aleph_0$-prebalancedness
(see \cite{Fu-new}). We are rephrasing the original definition in a form
which is clearly equivalent to the original definition.
\begin{definition} Let $G$ be a pure subgroup of the group $H$. $G$ is said to
be $\aleph_0$-prebalanced in $H$ if for every element $h\in H-G$ there
are countably many elements $g_0,g_1,\ldots$ of $G$ such that for
every element $g$ of $G$ the type (in $H$) of $h-g$ is bounded by the
the union of finitely many  types of the form $lh-g_i$ for some natural
number $l$. More
 explicitly for
some $n,l\in\omega$
$$\boldt(h-g)\leq\boldt(lh-g_0)\cup\ldots\cup\boldt(lh-g_n).$$
\end{definition}
Also the group $G$ is said to admit an $\aleph_0$-prebalanced chain if $G$
can be represented as a continuous increasing union of pure
$\aleph_0$-prebalanced subgroups where at the successor stages the factors
are of rank 1.

We shall use the following fundamental result of Fuchs (\cite{Fu-new}):
\begin{theorem}\label{Fu-main} A torsion-free group $G$ admits an
$\aleph_0$-prebalanced
chain if and only if in its balanced projective resolution
$$0\rightarrow B\rightarrow C\rightarrow G\rightarrow 0$$
(where $C$ is completely decomposable) $B$ is a $B_2$-group. Moreover, if CH
holds, then this condition is equivalent to $\bext=0$ for all torsion groups
$T$.
\end{theorem}

The main result of Section 3 will be
\begin{theorem}\label{main2} If (*) holds, then there is a
 torsion-free group $G$ of
cardinality $\aleph_{\omega+1}$ which does not admit an
$\aleph_0$-prebalanced chain.
\end{theorem}

Using theorem \ref{Fu-main} we get
\begin{cor}If (*) holds, then there is a group $G$ of cardinality
$\aleph_{\omega+1}$ such that $\bext\not=0$ for some torsion group $T$.
\end{cor}
By  using the balanced projective resolution of $G$ we also get
\begin{cor} If (*) holds, then there is a balanced subgroup of a completely
decomposable group of cardinality $\aleph_{\omega+1}$ which is not
a $B_2$-group.
\end{cor}
\section{The Consistency of (*)}
In this section we shall prove Theorem \ref{main1}. We assume familiarity
with some basic large cardinals notions like supercompact cardinals and
some basic forcing techniques. We start from a ground model $V$ having
a supercompact cardinal $\kappa$. We can assume without loss of generality
that $V$ satisfies G.C.H.
We let $\mu=\kappa^{+\omega}$ and $\lambda=\mu^+=\kappa^{+\omega+1}$.
In our final model $\mu$ will be $\aleph_\omega$ and
 $\lambda$ will be $\aleph_{\omega+1}$.
It follows from the results of  Menas in \cite{Menas} that
there is a normal ultrafilter $U$ on $P_\kappa(\lambda)$ such that for
some set $A\in U$ the map $P\rightarrow \sup(P)$ on $A$ is one-to-one.
(Recall that $P_\kappa(\lambda)$ is the set of all subsets of $\lambda$
of cardinality less than $\kappa$).
Fix such $U$ and $A$. Also fix an enumeration $\langle g_\alpha\mid
\alpha<\lambda\rangle$ of all the $\omega$-sequences in $\mu$. Standard
facts about normal ultrafilters on $P_\kappa(\lambda)$ imply that the set of
all $P\in P_\kappa(\lambda)$ satisfying the following properties is in $U$:
\begin{enumerate}
\item The order type of $P\cap\mu$ is a singular cardinal of cofinality
$\omega$ such that the order type of $P$ is its successor.
\item For $\alpha\in\lambda$ the range of $g_\alpha$ is a subset of
$P\cap\mu$ if and only if $\alpha\in P$.
\end{enumerate}
Hence we can assume without loss of generality that every $P\in A$ satisfies
all the above properties. Again standard arguments show that the
set $T=\{\sup(P)\mid P\in A\}$ is a stationary subset of $\lambda$. For
$\alpha\in T$, let $P_\alpha$ be the unique $P\in A$ such that
$\sup(P)=\alpha$. Note that for $P\in A$ and $Q\subseteq P$ we have that if
$Q$ is cofinal in $\sup(P)$, then the order type of
$Q^\ast=\cup\{range(g_\alpha)\mid\alpha\in Q\}$ is
the same as the order type of $P\cap\mu$. This holds since
otherwise $Q^\ast$ has cardinality smaller than $\delta=$the order
 type of $P\cap\mu$.
Hence, by our G.C.H. assumption, we have less than $\delta$ $\alpha$'s such
that the range of $g_\alpha$ is in $Q^\ast$, hence less than the order type
of $P$, which is a regular cardinal. Therefore $Q$ must be bounded in $P$.

For $\alpha\in T$ the map $\alpha\rightarrow\mbox{the order type of\
}P_\alpha\cap\mu$ maps $T$ into $\kappa$. Hence it is fixed on some subset
$S$ which is still stationary in $\lambda$. Let $\delta$ be the fixed value
of this map on $S$. Note that for $\alpha\in S$ the order type of $P_\alpha$
is $\delta^+$.
\begin{claim}\label{stype} Let $\alpha$ and $\beta$ be two different members
of
$S$. Then
$P_\alpha\cap  P_\beta\cap\mu$ has order type less than $\delta$.
\end{claim}
\begin{pf}
Let $X=P_\alpha\cap P_\beta\cap\mu$. Note that if $g$ is an $\omega$-sequence
 from
$X$, then $g=g_\rho$ for some $\rho\in P_\alpha\cap P_\beta$. If $X$ has order
type $\delta$, then (using the fact that $\delta$ is a singular cardinal
of cofinality $\omega$) we have $\delta^+$ $\omega$-sequences from $X$, so
that $P_\alpha\cap P_\beta$ must have order type which is at least
$\delta^+$. Since the order type of both $P_\alpha$ and $P_\beta$ is
$\delta^+$, $P_\alpha$ and $P_\beta$ must have the same sup. This is a
contradiction.
\end{pf}

The model which will  witness  (*) will be obtained from V by
collapsing $\delta$ to be countable, followed by the collpasing all the
cardinals between $\delta^{++}$ and $\kappa$ to have cardinality
$\delta^{++}$. Denote the resulting model by $V_1$. Note since $V$ satisfies
G.C.H. then the resulting
model satisfies G.C.H.\@ Also $\delta$ is of course countable, $\delta^+$
is $\aleph_1$ , $\mu$ is $\aleph_\omega$ and $\lambda$ is
$\aleph_{\omega+1}$.
Since the cardinality of the forcing notion is $\kappa<\lambda$, $S$ is
still a stationary subset of $\lambda$. Note that now we have for every
$\alpha\in S$ that the cofinality of $\alpha$ is $\aleph_1$. In order to
verify (*) in the
resulting model we fix an enumeration $\langle f_\gamma\mid
\gamma<\lambda\rangle$ of all the $\omega$-sequences from $\aleph_\omega=\mu$.
And as in the previous section let $F_\gamma$ be the range of $f_\gamma$.
(Note that in $V_1$ there are new $\omega$-sequences so that the enumeration
$\langle g_\gamma\mid \gamma<\lambda\rangle$ we had in $V$ enumerates only
a subset of the set of all $\omega$-sequences). For $\gamma<\lambda$ let
$\eta(\gamma)$ be the unique $\eta$ such that $g_\gamma=f_\eta$. Without
loss of generality (by reducing $S$ to a subset which is still
stationary in $\lambda$) we can assume that for $\alpha\in S$ if
$\gamma<\alpha$,
then
$\eta(\gamma)<\alpha$. We can also assume without loss of generality that
for $\alpha\in S$, $Q_\alpha=\{\eta(\gamma)\mid\gamma\in P_\alpha\}$ is
cofinal in $\alpha$. This follows since the set $\{\alpha\in S\mid
Q_\alpha\mbox{\ is bounded in\ } \alpha\}$ is not stationary. So for each
$\alpha\in S$
pick $C_\alpha$ which is cofinal in $Q_\alpha$ and has order type
$\aleph_1=\delta^+$. We claim that $S$, $\delta$ and $\langle
C_\alpha\mid\alpha\in S\rangle$ are witnesses to the truth of (*) in
$V_1$. As in the introduction we put
$$E_\alpha=\bigcup_{\gamma\in C_\alpha}F_\gamma.$$ Since we clearly have
G.C.H. in $V_1$, since $S$ is stationary and
since $\delta$ is an indecomposable ordinal (it is a cardinal in $V$!),
we are left with verifying the following claim:
\begin{claim}\label{chang} In $V_1$
\begin{description}
\item[A] For $\alpha\not=\beta\in S$ $E_\alpha\cap E_\beta$ has order type
less than $\delta$.
\item[B] If $D\subseteq C_\alpha$ is cofinal in $\alpha$, then
$\cup\{F_\gamma\mid \gamma\in D\}$ has order type $\delta$.
\end{description}
\end{claim}
\begin{pf} Clause A follows immediately from the fact that
for $\alpha\in S$, $E_\alpha\subseteq P_\alpha\cap\mu$, hence
$E_\alpha\cap E_\beta\subseteq P_\alpha\cap P_\beta\cap\mu$
and the last set has order type less than $\delta$ if $\alpha\not=\beta$.

For proving B note that if $D\subseteq C_\alpha$ is cofinal in $\alpha$,
then the set $F=\{\gamma\mid \eta(\gamma)\in D\}$ is a subset of $P_\alpha$
of cardinality $\aleph_1=\delta^+$. Our forcing is an iteration of two forcing
notions where the first is of cardinality (in $V$) $\delta$ and the second
is $\delta^{++}$ closed, hence it introduces no new sets of ordinals of
order type $\delta^+$. So $F$ contains a subset $Q\in V$ of cardinality
$\delta^+$. $Q$ must be cofinal in $P_\alpha$ since $P_\alpha$ has order type
$\delta^+$, so by a previous remark $\cup\{range(g_\gamma)\mid \gamma\in Q\}$
has order type $\delta$. But this last set is clearly a subset of
$\cup\{F_\rho\mid \rho\in D\}$, so this set clearly has order type at
least $\delta$. It can not have order type greater than $\delta$ since it is
a subset of $P_\alpha\cap\mu$.
\end{pf}
\section{A group which does not admit an $\aleph_0$-prebalanced chain}
 In this section we prove Theorem \ref{main2}. So we assume (*). Fix the
 enumeration $\langle f_\alpha\mid \alpha<\aleph_{\omega+1}\rangle$
 of the $\omega$-sequences from $\aleph_\omega$. Let $F_\alpha$ be the range
 of $f_\alpha$. Also fix
 the stationary subset $S$ of $\aleph_{\omega+1}$, the countable ordinal
 $\delta$ and for $\beta\in S$ a set $C_\beta$ cofinal in $\beta$, which
 witness the truth of (*). As in the statement of (*) (for $\beta\in S$)
 let
 $$E_\beta=\bigcup_{\alpha\in C_\beta}F_\alpha.$$ We know that the order
 type of $E_\beta$ is $\delta$. Since $\delta\times\omega$ is countable
 we can assign to every pair $\mu<\delta,n<\omega$ a unique prime number
 $\pnmu$.

 We are ready to define the group $G$ that will not admit a chain of
 $\aleph_0$-prebalanced subgroups. For each $\alpha<\aleph_{\omega+1}$
   and $\beta\in S$
  fix distinct symbols $x_\alpha$ and $y_\beta$.
 The group $G$ is a subgroup of
 $$\sum_{\alpha<\aleph_{\omega+1}}\!\!\!\!\oplus \boldQ
  x_\alpha\oplus\sum_{\beta\in S}\!\!\oplus
 \boldQ y_\beta.$$ $G$ is generated by $x_\alpha$ for
  $\alpha<\aleph_{\omega+1}$, by
  $y_\beta$ for $\beta\in S$ and by
   $\frac{1}{\pnmu}(y_\beta-x_\alpha)$ provided
  $\alpha$ is in $C_\beta$ and the $f_\alpha(n)$ is the $\mu$-th
  member of $E_\beta$. For $\delta<\aleph_{\omega+1}$ let $G_\delta$ be the
  subgroup of $G$ generated by $x_\alpha$,$y_\gamma$ and
  $\frac{1}{\pnmu}(y_\gamma-x_\alpha)$ where $\alpha$ and $\gamma$ are less
  than $\delta$. The sequence $\langle G_\delta\mid
  \delta<\aleph_{\omega+1}\rangle$
  is a filtration of $G$ into a continuous chain of smaller cardinality.
  If $G$ allows an $\aleph_0$-prebalanced chain, then by standard
  arguments, the set of $\delta<\aleph_{\omega+1}$ such that $G_\delta$
  appears in the $\aleph_0$-prebalanced chain contains a closed unbounded
  subset of $\aleph_{\omega+1}$. This will imply, since $S$ is stationary
  in $\aleph_{\omega+1}$, that for some $\beta\in S$, $G_\beta$ is $\aleph_0$
  prebalanced in $G$. The fact that we get a contradiction and that $G$
  does not allow an $\aleph_0$-prebalanced chain follows from:
  \begin{claim}\label{mclaim} For $\beta\in S$, $G_\beta$ is not an
$\aleph_0$-prebalanced
  subgroup of $G$.
  \end{claim}
  \begin{pf}
  Assume that for some fixed $\beta\in S$, $G_\beta$ is
  $\aleph_0$-prebalanced in $G$. We apply the definition of
  $\aleph_0$-prebalancedness for $y_\beta$ and get a sequence of elements
  $z_n\in G_\beta$ such that for every element $z$ of $G_\beta$ there are $e$
  and $l$
  such that
 $$\boldt(y_\beta-z)\leq\boldt(ly_\beta-z_0)\cup\ldots\cup%
 \boldt(ly_\beta-z_e).$$
 $C_\beta$ has order type $\aleph_1$ and hence for some fixed $e$ and $l$
   we get that
   the set
  \begin{equation}\label{type}
  D=\{\alpha\in C_\beta\mid \boldt(y_\beta-x_\alpha)\leq
  \boldt(ly_\beta-z_0)\cup\ldots\cup\boldt(ly_\beta-z_e)\}
  \end{equation}
  is unbounded in $C_\beta$. It means that for $\alpha\in D$ there is a
  natural number $d_\alpha$ such that if $p$ is a prime number greater than
  $d_\alpha$ and $p$ divides $y_\beta-x_\alpha$, then $p$ divides
  $ly_\beta-z_i$ for some $0\leq i\leq e$. Without loss of generality we can
  assume that for $\alpha\in D$, $d_\alpha$ is some fixed natural number $d$.
  Let $D^\ast=\cup_{\gamma\in D}F_\gamma$. We
  know that $D^\ast\subseteq E_\beta$ and that the order type of $D^\ast$ is
  $\delta$.
  We need the following lemma.
  \begin{lemma}\label{mlemma} Let $z$ be a member of $G_\beta$ with
  $$z=\sum_{i=1}^{k}r_ix_{\alpha_i}+\sum_{j=1}^{g}s_jy_{\beta_j},$$
  where $r_i,s_j\in\boldQ$ and $\alpha_i,\beta_j<\beta$ for
  $1\leq i\leq k,1\leq j\leq g$. Assume also that $ly_\beta-z$ is divisible
  (in $G$) by $\pnmu$ where $\pnmu>l$. Then either for some $1\leq j\leq g$,
   the
  $\mu$-th member of $E_\beta$ is the same as the $\mu$-th member
  of $E_{\beta_j}$ or for some $1\leq i\leq k$, the $\mu$-th member
  of $E_\beta$ is in $F_{\alpha_i}$.
  \end{lemma}
  \begin{pf} By assumption $ly_\beta-z$ is divisible by $p=\pnmu$ in $G$.
  Hence \begin{equation}\label{pdiv}
  ly_\beta-z=p(\sum_{m=1}^{f}r_mx_{\gamma_m}+\sum_{t=1}^{u}s_ty_{\eta_t}+
  \sum_{q=1}^{v}\frac{w_q}{p_q}(y_{\nu_q}-x_{\xi_q})).
  \end{equation}
  where the $r_m$'s , the $s_t$'s and the $w_q$'s are integers.

  Let us define a (bipartite) graph $P$, whose nodes are all the symbols
  ($x$'s and $y$'s) appearing in equation \ref{pdiv}, where $y_{\rho}$ is
  connected by an edge to $x_{\zeta}$ iff for some $1\leq q \leq v,$
  $\rho=\nu_q$, $\zeta=\xi_q$ and $p_q=p$. Let $W$ be the connected component
  of $y_\beta$ in $P$ and let $a\in\boldQ$ be the sum of all the
   coefficients
  in the right  side of equation \ref{pdiv} of symbols in $W$.
   $a$  is easily
  seen to be a member of $p\boldQ_p$, where $\boldQ_p$ is
  the ring of rationals whose
  denominators are prime to $p$.This is true because
   the only summands on the right side of
  \ref{pdiv}, that can possibly add to $a$ a rational number which is not in
  $p\boldQ_p$, is of the form $\frac{w_q}{p_q}(y_{\nu_q}-x_{\xi_q})$ where
  $p_q=p$. But in this case  $y_{\nu_q}$ and $x_{\xi_q}$ are connected by an
  edge of $P$, so they are both in $W$ or both outside of $W$. In both cases
  the contribution of this summand to $a$ is $0$.

  We use the fact that the sum of the coefficients of symbols in $W$
  must be the same for the left side and the right side of \ref{pdiv}.
    Of course $y_\beta\in W$ and its coefficient in equation~\ref{pdiv} is $l$
which is
  not in $p\boldQ_p$, so there must be a symbol in $W$ appearing in the
  representation of $z$, so that either $x_{\alpha_i}\in W$ for some
  $1\leq i\leq k$, or $y_{\beta_j}\in W$ for some $1\leq j\leq g$. Our
  lemma will be verified if we prove
  \begin{claim}\label{pclaim} \begin{enumerate}
  \item If $y_\eta\in W$, then the $\mu$-th member of $E_\eta$ is the
  same as the $\mu$-th member of $E_\beta$.
  \item If $x_\gamma\in W$, then $f_\gamma(n)$ is the $\mu$-th member of
  $E_\beta$.
  \end{enumerate}
  \end{claim}
  \begin{pf}
   The proof is by induction on the length of the path in $P$ leading from
   $y_\beta$ to the symbol $y_\eta$ and $x_\gamma$ respectively. If this
   length is $0$, we are in the case where the symbol is $y_\eta=y_\beta$, and
   the claim is obvious. For the induction step, in the first case we are
   given $y_\eta$. Let $x_\gamma$ be the element preceding $y_\eta$ in the
   path leading from $y_\beta$ to $y_\eta$. By the induction assumption
   $f_\gamma(n)$ is the $\mu$-th member of $E_\beta$. $x_\gamma$ and $y_\beta$
   are connected by an edge of $P$, so that
   $\frac{1}{\pnmu}(y_\eta-x_\gamma)$ is one of the generators of $G$.
   Hence $\gamma\in C_\eta$ and $f_\gamma(n)$ is the $\mu$-th member
   of $E_\eta$, and the claim is verified in this case. The other case (the
   $x_\gamma$ case) is argued similary where $y_\eta$ is now the element in
   the path preceding $x_\gamma$.
  \end{pf}
  \end{pf}

 For $z\in G_\beta$ let $S(z)$ be the set of all elements $\gamma$ of
 $E_\beta$ such that for some $\mu<\delta$ and $n\in\omega$,
 $\gamma$ is the $\mu$-th member of $E_\beta$ and $ly_\beta-z$ is divisible
 in $G$ by $\pnmu$ where $\pnmu>l$.
 It follows from lemma~\ref{mlemma} that for $z\in G_\beta$, $S(z)$ is
 included in
  a finite
 union of singletons and of sets of the form $E_\eta\cap E_\beta$ for
 $\eta<\beta$. So $S(z)$ is a finite union of sets of order type less than
 $\delta$. $\delta$ is an indecomposable ordinal, so for
   $z\in G_\beta$ the order type of $S(z)$ is less than $\delta$.
   By definition of $D$, every element of $D^\ast$, except possibly finitely
   many, is in $\cup_{0\leq i\leq e}S(z_i)$. This is because there are only
   finitely many members of $E_\beta$ such that if $\gamma$ is the
   $\mu$-th member of $E_\beta$, then $\pnmu\leq \max(d,l)$ for
    some $n$. So if $\gamma\in
   D^\ast$ is not one of these finitely many elements, say $\gamma$ is the
   $\mu$-th member of $E_\beta$, then $\pnmu> \max(d,l)$. Now $
   \gamma=f_\alpha(n)$
   for some $\alpha\in D$ and a natural number $n$, and hence $\pnmu$ divides
   $y_\beta-x_\alpha$, which implies by equation~\ref{type} and the
   definition of $d$ that $\pnmu$ divides $ly_\beta-z_i$ for some $1\leq
   i\leq e$. We got that $D^\ast$ is a finite union of sets of order type
   less that $\delta$, and hence $D^\ast$ has order type less than $\delta$.
   We got a contradiction.
   \end{pf}
   \bibliography{groups}

\ifx\undefined\bysame
\newcommand{\bysame}{\leavevmode\hbox to3em{\hrulefill}\,}
\fi
\begin{thebibliography}{10}

\bibitem{Al-Hi}
U.~Albrecht and P.~Hill, {\em Butler groups of infinite rank}, Czech. Math. J.
  {\bf 37} (1987), 293--309.

\bibitem{Bi-Sa}
L.~Bican and L.~Salce, {\em Infinite rank {Butler } groups}, Lecture Notes In
  Mathematics, vol. 1006, Springer, 1983.

\bibitem{DHR}
M.~Dugas, P.~Hill, and K.~M. Rangaswamy, {\em Infinite rank {Butler} groups
  {II}}, Trans. Amer. Math. Soc. {\bf 320} (1990), 643--664.

\bibitem{DR}
M.~Dugas and K.~M. Rangaswamy, {\em Infinite rank {Butler} groups}, Trans.
  Amer. Math. Soc. {\bf 305} (1988), 129--142.

\bibitem{Du-Th}
M.~Dugas and R.~Thom\'e, {\em The functor {Bext} under the negation of {CH}},
  Forum Math. {\bf 3} (1991), 23--33.

\bibitem{Ek-Me}
P.C. Eklof and A.H.Mekler, {\em Almost free modules}, North-Holland, 1990.

\bibitem{For-Ma}
M.~Foreman and M.~Magidor, {\em A version of weak {$\Box$}}, to appear.

\bibitem{Fu-new}
L.~Fuchs, {\em A survey of {Butler} groups of infinite rank}, this volume.

\bibitem{Fu.book}
\bysame, {\em Infinite abelian groups}, vol.~I, Academic Press, New York and
  London, 1970.

\bibitem{Fu-Ma}
L.~Fuchs and M.~Magidor, {\em Butler groups of arbitrary cardinality}, Israel
  Journal of Math. {\bf 84} (1993), 239--263.

\bibitem{HJSh249}
A.~Hajnal, I.~Juhasz, and S.~Shelah, {\em Splitting strongly almost disjoint
  families}, Trans. of the A.M.S. {\bf 295} (1986), 369--387.

\bibitem{Menas}
T.K.Menas, {\em A combinatorial property of {$P_\kappa\lambda$}}, J. of
  Symbolic Logic {\bf 41} (1976), 225--234.

\end{thebibliography}
   \end{document}